\documentclass{article}

\usepackage{arxiv}

\usepackage[utf8]{inputenc} 
\usepackage[T1]{fontenc}    
\usepackage{hyperref}       
\usepackage{url}            
\usepackage{booktabs}       
\usepackage{amsfonts}       
\usepackage{nicefrac}       
\usepackage{microtype}      
\usepackage{lipsum}		
\usepackage{graphicx}
\usepackage{natbib}
\usepackage{doi}
\usepackage{amsmath}

\newtheorem{Theorem}{Theorem}
\newtheorem{Lemma}[Theorem]{Lemma}
\newtheorem{Definition}[Theorem]{Definition}
\newtheorem{Corollary}[Theorem]{Corollary}

\title{On the cycle maximum of birth-death processes and networks of queues}

\author{ \href{https://orcid.org/0000-0002-1046-2044}{\includegraphics[scale=0.06]{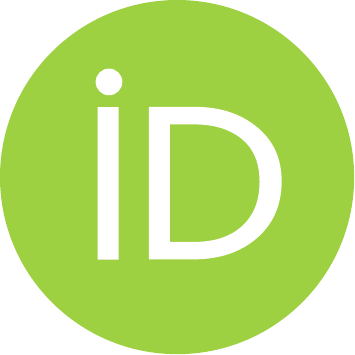}\hspace{1mm}Richard J. Boucherie}
\\
	Stochastic Operations Research \\ Department of Applied Mathematics\\
	University of Twente \\
	\texttt{r.j.boucherie@utwente.nl} \\
}

\renewcommand{\headeright}{Technical Report}
\renewcommand{\undertitle}{Technical Report}
\renewcommand{\shorttitle}{On the cycle maximum of birth-death processes and networks of queues}

\hypersetup{
pdftitle={A template for the arxiv style},
pdfsubject={q-bio.NC, q-bio.QM},
pdfauthor={David S.~Hippocampus, Elias D.~Striatum},
pdfkeywords={First keyword, Second keyword, More},
}

\begin{document}
\maketitle

\begin{abstract}
 This paper considers the cycle maximum in birth-death processes as a stepping stone to characterisation of the cycle maximum in single queues and open Kelly-Whittle networks of queues. For positive recurrent birth-death processes we show that the sequence of sample maxima is stochastically compact. For transient birth-death processes we show that the sequence of sample maxima conditioned on the maximum being finite is stochastically compact. 
 
 We show that the Markov chain recording the total number of customers in a Kelly-Whittle network is a birth-death process with birth and death rates determined by the normalising constants in a suitably defined sequence of closed networks. Explicit or asymptotic expressions for these normalising constants  allow asymptotic evaluation of the birth and death rates, which, in turn, allows characterisation of the cycle maximum in a single busy cycle, and convergence of the sequence of sample maxima for Kelly-Whittle networks of queues. 
\end{abstract}


\section{Introduction}
The maximum number of customers simultaneously present in a queueing system is an important performance measure that is directly related to the cycle maximum during a busy cycle (\cite{Asmussen1998}). This paper considers the cycle maximum in birth-death processes as a stepping stone to characterisation of the cycle maximum in single queues and open Kelly-Whittle networks of queues (\cite{BoucherievanDijk1991,HendersonTaylor1990}). To this end, it is observed that Norton's theorem (\cite{BoucherieNorton}) implies that the Markov chain recording the total number of customers in the network is a birth-death process with birth and death rates determined by the normalising constants in a suitably defined sequence of closed networks. Explicit or asymptotic expressions for these normalising constants  (\cite{Harrison1985,GeorgeXiaSquillante2012}) allow asymptotic evaluation of the birth and death rates, which, in turn, allows characterisation of the cycle maximum in a single busy cycle, ${\mathbf Y}$, and convergence of the sequence of partial maxima, ${\mathbf Y}^{(k)}$, over $k$ busy cycles. 

Building on the explicit expression for the cycle maximum in a birth-death process (\cite{Asmussen1998}), we observe that the dual relation between the discrete failure rate $
\mathbb{P}(\mathbf{Y} = n | \mathbf{Y} \geq n)$ and the blocking probability $
\mathbb{P}(\mathbf{X} = n | \mathbf{X} \leq n)$, with $\mathbf X$ the limiting random variable for state of the process, that was first observed  for the Erlang loss queue in \cite{Cohen1971}, extends to general birth-death processes. This dual relation is of a form similar to the dual or inverse relation between the $M/G/1$ and $GI/M/1$ queues obtained by interchanging the service and inter-arrival time distributions, see \cite{Takacs1962,NiuCooper1989,Kimura1993,BekkerZwart2005}. 

The limiting behaviour of the sequence of sample maxima ${\mathbf Y}^{(k)}$ may be characterised by a sequence of constants $\{a_k>0,~b_k\}_{k\geq 1}$, $a_k$ bounded, $b_k\rightarrow \infty$ as $k\rightarrow \infty$. For continuous valued $\mathbf Y$ such as for the waiting time or workload these constants are available in literature, see e.g. \cite{BoxmaPerry2009}. As $\mathbf Y$ for the birth-death process is discrete valued, $\lim_{k \rightarrow \infty} \mathbb{P}(\mathbf{Y}^{(k)} \leq a_k x + b_k)$ does not exist, see \cite{Serfozo1988}. Bounds for the $\liminf$ and $\limsup$ may be obtained following the approach in \cite{Anderson1970,Cohen1982,Vervaat1973}. For positive recurrent birth-death processes, we show that $\mathbf{Y}$ is stochastically compact (\cite{DeHaanResnick1984}), in particular, we will show that the normalised sequence of sample maxima $(\mathbf{Y}^{(k)} - b_k)/ a_k$ has a subsequence whose distributions converge weakly to the Gumbel distribution. For transient birth-death processes, we show that $\mathbf{Y} | \mathbf{Y}<\infty$  is stochastically compact. The relation between  $\mathbf{Y}$ for recurrent processes and  $\mathbf{Y} | \mathbf{Y}<\infty$  for transient processes is similar to the relation first observed in \cite{Cohen1971} between a process and its dual, as mentioned above. 

Section \ref{sec:2.1} introduces notation for birth-death processes, and observes the dual relation for the cycle maximum. Section \ref{sec:3} considers weak and almost sure convergence of sample maxima. Section \ref{sec:4} considers the cycle maximum for single queues and a Kelly-Whittle network of queues.

\section{The cycle maximum for a single cycle}\label{sec:2.1}
\subsection{Birth-death process}
Consider a birth-death process $(X(t),~t \geq 0)$ at state space $S= \{0,1,2,\ldots\}$ with birth rates $\lambda(n)>0$,  $n\geq 0$ and death rates $\mu(n)>0$, $n\geq 1$. In agreement with common notation for queues and queueing networks, see \cite{BoucherievanDijk1991,HendersonTaylor1990}, we will express the birth and death rates as 
\begin{equation}
\lambda(n)= \lambda \frac{\psi(n)}{\phi(n)}, \quad n \geq 0, 
\quad \mu(n)= \mu \frac{\psi(n-1)}{\phi(n)}, \quad n>1,
\label{eq:rates}
\end{equation}
with $\psi(n)>0$, $\phi(n)>0$, $n=0,1,2,\ldots$. Let $\rho=\lambda / \mu$.

We assume that $\sum_{n\geq 1} \phi(n) / (\psi(n)+\psi(n-1))=\infty$, which is a sufficient condition for regularity (\cite{Foster1953}). $(X(t),~t \geq 0)$  is positive recurrent if and only if (iff) $B_{\phi}^{-1}=\sum_{n\in S} \phi(n)\rho^n<\infty$. Then the limiting distribution $\pi_{\phi}(n)=\lim_{t\rightarrow \infty}\mathbb{P}(X(t)=n)$ exists, and $\pi_{\phi}(n)=B_{\phi}\phi(n)\rho^n$, $n \in S$. Let $\mathbf X$ denote the limiting random variable. From Foster's criterion (\cite{Foster1953}), $(X(t),~t \geq 0)$  is null recurrent iff $B_{\phi}^{-1}=\infty$ and $B_{*}^{-1}=\sum_{n\in S}(\psi(n)\rho^n)^{-1}<\infty$
 and transient iff $B_{*}^{-1}=\infty$. Therefore, if  $B_{\psi}^{-1}=\sum_{n\in S}\psi(n)\rho^n<\infty$, then $(X(t),~t \geq 0)$  recurrent. Furthermore, if $B_{\psi}^{-1}<\infty$   the limiting Palm distribution at births $\pi_{\psi}(n)=\lim_{t\rightarrow \infty}\mathbb{P}_0(X(t)=n)=\lim_{t\rightarrow \infty} \lim_{\Delta t \downarrow 0} \mathbb{P}(X(t)=n | \mbox{ birth in } (t,t+\Delta t) )
 $ exists, and $\pi_{\psi}(n) =B_{\psi}\psi(n)\rho^n$, $n \in S$, 
see \cite{Brumelle1978}.

Observe that  $(X(t),~t \geq 0)$ is regenerative  with renewal epochs $0\leq T_0<T_1<T_2< \cdots$ each time $(X(t),~t \geq 0)$ has a transition from state 0 to state 1. The generic cycle is the cycle starting at $T_0$. Let $\mathbf{Y}_k = \max_{T_{k-1}\leq t <T_k} X(t)$ the maximum over the $k$-th cycle. Then $\mathbf{Y}_1, \mathbf{Y}_2,\ldots $ are i.i.d.\ and distributed as $\mathbf{Y} = \mathbf{Y}_1$, with distribution, see \cite{Asmussen1998},
\begin{equation}
\mathbb{P}(\mathbf{Y} \leq n)= 1-\left[ \sum_{i=0}^n \frac{1}{\psi(i)\rho^i}\right]^{-1}, \quad n=1,2,\ldots,
\label{eq:Y}
\end{equation}
which is non-degenerate iff $B_*^{-1}=\infty$, i.e., iff $(X(t),~t \geq 0)$ is recurrent.

For the $M(\lambda)/M(\mu)/1$-queue, the single server queue with Poisson arrivals with rate $\lambda$ and exponential service requirements with rate $\mu$,  $\psi(n)=\phi(n)=1$. For the $M(\lambda)/M(\mu)/ \infty$-queue, the infinite-server queue,  $\psi(n)=\phi(n)= 1  / {n!}$. See Section \ref{sec:4} for further details.

The result for the cycle maximum readily extends to a birth-death process with finite state space $S=\{0,1,\ldots,c\}$, which may be achieved by setting $\psi( c )=0$. For example, consider the $M(\lambda)/M(\mu)/1/c$-queue, the single server queue finite waiting room of size $c-1$, for which the distribution of the cycle maximum is obtained in \cite{HanbaliBoxma2010}, or the  $M(\lambda)/M(\mu)/c/c$-queue, the Erlang loss queue.  The probability $\mathbb{P}(\mathbf{Y}\leq c)$ must now be interpreted as the probability that during a cycle no customers are blocked.

\subsection{A duality relation for the cycle maximum  for queues in equilibrium}

From (\ref{eq:Y}) we obtain the conditional probability that the cycle maximum is $n$ given that the cycle maximum is at least $n$:
\begin{equation}
\mathbb{P}(\mathbf{Y} = n | \mathbf{Y} \geq n)= \frac{1}{\psi(n)\rho^n} \left[ \sum_{i=0}^n \frac{1}{\psi(i)\rho^i}\right]^{-1}. \label{eq:Y|Y}
\end{equation}
From the Palm distribution $\mathbb{P}_0(\mathbf{X} = n)$ we obtain the probability that the state at a birth is $n$ given that it is at most $n$:
\begin{equation}
\mathbb{P}_0(\mathbf{X} = n | \mathbf{X} \leq n)= {\psi(n)\rho^n} \left[ \sum_{i=0}^n {\psi(i)\rho^i}\right]^{-1}.
\label{eq:X|X}
\end{equation}
The relation between (\ref{eq:Y|Y}) and (\ref{eq:X|X})  was first observed in \cite{Cohen1971} for the Erlang loss queue, where $\mathbb{P}_0(\mathbf{X} = c | \mathbf{X} \leq c)$ is the Erlang loss formula, i.e., the probability that an arriving customer meets $c$ customers and is blocked, and $\mathbb{P}(\mathbf{Y} = c | \mathbf{Y} \geq c)$ is the probability that no customer is blocked during a cycle, whenever during this cycle all servers are simultaneously busy at least once. 

For the $M(\lambda)/M(\mu)/1$-queue,  with stability condition $\lambda< \mu$, the relation  between (\ref{eq:Y|Y}) and (\ref{eq:X|X}) is more prominent. We find that   
\[
\mathbb{P}_{M(\mu)/M(\lambda)/1}(\mathbf{Y} = n | \mathbf{Y} \geq n)= \mathbb{P}_{M(\lambda)/M(\mu)/1}(\mathbf{X} = n | \mathbf{X} \leq n),
\]
where the non-stable $M(\mu)/M(\lambda)/1$-queue is the dual of the stable $M(\lambda)/M(\mu)/1$-queue, i.e., the queue obtained by interchanging the inter-arrival and service time distributions (\cite{Kimura1993}).

\section{Limit theorems}\label{sec:3}
This section first considers the tail behaviour of $\mathbb{P}(\mathbf{Y} \leq n)$, and subsequently analyses the normalised sample maxima $\mathbb{P}((\mathbf{Y}^{(k)}-b_k)/a_k \leq x)$ for $\mathbf{Y}^{(k)}=\max\{\mathbf{Y}_1,\ldots,\mathbf{Y}_k\}$, and suitable constants $\{a_k>0,~b_k\}_{k\geq 1}$. We will characterise weak and almost sure convergence of the sample maxima for recurrent processes.

\subsection{Tail behaviour}
Recurrence and transience of $(X(t),~t \geq 0)$ is characterised by $\psi$, see Section \ref{sec:2.1}.  The ratio test for convergence of the series $B_{\psi}^{-1}$ and $B_*^{-1}$ suggest the following notation to characterise convergence: 
\begin{eqnarray*} 
\overline{\beta} &=& \limsup _{n\rightarrow \infty} \psi(n+1)/ \psi(n) , \\
\underline{\beta} &=& \liminf _{n\rightarrow \infty} \psi(n+1)/ \psi(n) ,\\
{\beta} &=& \lim_{n\rightarrow \infty} \psi(n+1)/ \psi(n), \quad \mbox{if the limit exists.}
\end{eqnarray*} 
If $(X(t),~t \geq 0)$ is recurrent, then  $\lim_{n\rightarrow \infty} \mathbb{P}(\mathbf{Y} \leq n) = 1$. We have the following asymptotic result. 
\begin{Lemma}\label{Lem:1}
If $\overline{\beta}\rho <1$, then. 
\begin{equation}
1-\overline{\beta}\rho \leq \liminf_{n\rightarrow \infty} \frac{1-\mathbb{P}(\mathbf{Y} \leq n)}{\psi(n)\rho^n}\leq \limsup_{n\rightarrow \infty} \frac{1-\mathbb{P}(\mathbf{Y} \leq n)}{\psi(n)\rho^n}\leq 1-\underline{\beta} \rho.
\label{eq:3.4}
\end{equation}
If $\beta$ exists, and $\beta\rho<1$, then 
\begin{equation}
 \lim_{n\rightarrow \infty} \frac{1-\mathbb{P}(\mathbf{Y} \leq n)}{\psi(n)\rho^n}= 1-{\beta}\rho .
 \label{eq:3.5}
 \end{equation}
\end{Lemma}
{\bf Proof} The results follow since $T(n)= \frac{1-\mathbb{P}(\mathbf{Y} \leq n)}{\psi(n)\rho^n}$ satisfies 
\[
\frac{1}{T(n+1)}=\rho\frac{\psi(n+1)}{\psi(n)}\frac{1}{T(n)}+1,
\] and taking limits $n\rightarrow \infty$. \hfill $\Box$ \\

If $(X(t),~t \geq 0)$ is transient, then  $\lim_{n\rightarrow \infty} \mathbb{P}(\mathbf{Y} \leq n) = 1-B_*$. We have the following asymptotic result.
\begin{Lemma}
If $\beta$ exists, and $\beta\rho>1$, then 
\begin{equation}
 \lim_{n\rightarrow \infty} \psi(n)\rho^n [1-\mathbb{P}(\mathbf{Y} \leq n | \mathbf{Y}<\infty)]= [({\beta}\rho-1)B_*(B_*-1) ]^{-1}.
 \label{eq:3.8}
 \end{equation}
\end{Lemma}
{\bf Proof} The result follows, since $T(n)=\psi(n) \rho^n[1-\mathbb{P}(\mathbf{Y} \leq n | \mathbf{Y}<\infty)]$ satisfies 
\[
{T(n+1)}=\rho\frac{\psi(n+1)}{\psi(n)}\frac{\sum_{i=0}^{n} (\psi(i)\rho^i)^{-1}}{\sum_{i=0}^{n+1} (\psi(i)\rho^i)^{-1}}T(n)-\frac{1}{(B_*-1)\sum_{i=0}^{n+1} (\psi(i)\rho^i)^{-1}},
\]
 and taking limits $n\rightarrow \infty$. \hfill $\Box$ \\

From (\ref{eq:3.5}) and (\ref{eq:3.8}) we find, if $\beta$ exists, 
\begin{align}
& \lim_{n\rightarrow \infty} [1-\mathbb{P}(\mathbf{Y} \leq n)]/ [B_{\psi}\psi(n)\rho^n]= (1-{\beta}\rho)/B_{\psi} , & \beta\rho<1, \label{eq:Rem34-1}\\
&\lim_{n\rightarrow \infty} [1-\mathbb{P}(\mathbf{Y} \leq n | \mathbf{Y}<\infty)] / [1/(B_*\psi(n)\rho^n)]= 1/[(\beta\rho-1) \mathbb{P}(\mathbf{Y} <\infty) B_*], &  \beta \rho>1. \label{eq:Rem34-2}
\end{align}
If $\beta\rho<1$, then $B_{\psi}\psi(n)\rho^n$ is a non-degenerate distribution, and for $\beta\rho>1$, $1/(B_*\psi(n)\rho^n)$ is a non-degenerate distribution. This shows that the tail behaviour of $\mathbf{Y}$ conditioned on $\mathbf{Y}<\infty$ for transient birth-death processes is similar to that of $\mathbf{Y}$ for recurrent birth-death processes in the sense that $B_*\psi(n)\rho^n$ is the "dual" of $B_{\psi}\psi(n)\rho^n$ obtained by replacing the birth and death rates $\lambda(n)$ and $\mu(n)$ by $\lambda_*(n)=\mu \frac{\phi(n)}{\psi(n)}$ and $\mu_*(n)=\lambda \frac{\phi(n)}{\psi(n-1)}$.

If $\beta$ exists, and $\beta\rho=1$, the  tail behaviour of $\mathbb{P}(\mathbf{Y} \leq n)$ depends on tail behaviour of $\psi(n)$. 
\begin{Lemma}\label{Lem:3}
Assume  that $\lim_{n \rightarrow \infty} \psi(n)\rho^n/n^p = \alpha$ for some $\alpha$, $p$, with $0<\alpha<\infty$, and $-\infty < p < \infty$. Then, for all $p$, $-\infty < p < \infty$, there exist $\delta(n,p)$, $\gamma( p )$ such that
\[
\lim_{n\rightarrow \infty} \delta(n,p) [1-\mathbb{P}(\mathbf{Y} \leq n)] = \gamma (p ),
\]
where
\[
\delta(n,p) = \left\{ \begin{array}{ll} n^{1-p}, & \mbox{if }p<1,  \\ \log n , & \mbox{if } p=1, \\ 1, & \mbox{if } p>1, \end{array} \right. 
\quad \left\{ \begin{array}{ll} \alpha({1-p}), & \mbox{if }p<1,  \\ \alpha , & \mbox{if } p=1, \\ B_*, & \mbox{if } p>1. \end{array} \right. 
\]
\end{Lemma}
{\bf Proof} For all $\epsilon>0$, there exits an $n(\epsilon)$ such that $(\alpha-\epsilon)n^p < \psi(n)\rho^n < (\alpha+\epsilon)n^p$, so that for all $n>n(\epsilon)$
\[
\sum_{i=1}^{n\epsilon)} \frac{1}{\psi(i)\rho^i }+ \frac{1}{\alpha + \epsilon} \sum_{i=n(\epsilon)+1}^n \frac{1}{i^p} \leq [1-\mathbb{P}(\mathbf{Y} \leq n)] ^{-1} \leq \sum_{i=1}^{n\epsilon)} \frac{1}{\psi(i)\rho^i }+ \frac{1}{\alpha + \epsilon} \sum_{i=n(\epsilon)+1}^n \frac{1}{i^p}.
\]
 $\sum_{i}{i^{-p}}$ converges for $p>1$,  $\sum_{i=1}^n{i^{-p}}$ behaves like $\log n$ for $p=1$, and like $n^{1-p}$ for $p<1$, which yield the three cases and corresponding constants $\gamma( p)$. \hfill $\Box$ \\
 
Similar results may be obtained for other choices than $\psi(n) \rho^n\sim n^p$ for the tail behaviour of $\psi(n)$, such as  $\psi(n)\rho^n \sim n(\log n)^p$ and  $\psi(n)\rho^n \sim n \log n  (\log\log n)^p$.


\subsection{Weak convergence of sample maxima}
We will focus on the recurrent cases $\overline{\beta}\rho<1$, and if $\beta$ exists on the cases $\beta \rho<1$ and $ \beta \rho >1$, for which we show stochastic compactness. 

A necessary and sufficient condition for the sequence of sample maxima $\mathbf{Y}^{(k)}$ to possess normalisation constants leading to a non-degenerate limit  is that $\lim_{x \rightarrow x_F} \mathbb{P}(\mathbf{Y} > x) / \mathbb{P}(\mathbf{Y} \geq x) =1$ \cite[Theorem 1.7.13]{Leadbetteretal1983}. This condition is not satisfied for most discrete valued distributions such as the Poisson distribution \cite[Example 1.7.14]{Leadbetteretal1983}, and the geometric distribution  \cite[Example 1.7.15]{Leadbetteretal1983}, that are typically encountered as limiting distribution of queues. For such systems we may obtain the $\liminf$ and $\limsup$ for the distribution of the normalised sample maxima, similar to results presented in \cite{Anderson1970,Cohen1982}. 

\begin{Theorem}\label{Thm:3.5}
Assume that $0 \leq \overline{\beta} \rho <1$. Then there exists a continuous and decreasing function $f:[0,\infty) \rightarrow [0,\infty)$ such that for $-\infty < x < \infty$
\begin{align}
\liminf_{k \rightarrow \infty} \mathbb{P}(\mathbf{Y}^{(k)} \leq f^{-1}({\rm e}^{-x} / [(1-\underline{\beta}\rho)k]+1) \geq {\rm e}^{-{\rm e}^{-x}}, \label{eq:thm4-1}\\
\limsup_{k \rightarrow \infty} \mathbb{P}(\mathbf{Y}^{(k)} \leq f^{-1}({\rm e}^{-x} / [(1-\overline{\beta}\rho)k]) \leq {\rm e}^{-{\rm e}^{-x}}. \label{eq:thm4-2}
\end{align}
If in addition $\underline{\beta} >0$ then for $-\infty < x < \infty$
\begin{align}
\liminf_{k \rightarrow \infty} \mathbb{P}(\mathbf{Y}^{(k)} \leq f^{-1}({\rm e}^{-x} / [(1-\underline{\beta}\rho)k]) \geq {\rm e}^{-(\underline{\beta}\rho)^{-1}{\rm e}^{-x}}. \label{eq:thm4-3}
\end{align}
\end{Theorem}
{\bf Proof} 
Consider the continuous function $g:[0,\infty) \rightarrow [0,\infty)$ obtained from $\psi(n)\rho^n$ by linear interpolation. As $\liminf_{y\rightarrow\infty} g(y+1)/g(y)=\underline{\beta}$, $\limsup_{y\rightarrow\infty} g(y+1)/g(y)=\overline{\beta}$,  there exists an $y_0$ such that $g(y)$ is decreasing for $y>y_0$.  Set $f(y)=g(y)$, $y\geq y_0$, and $f(y)=y_0-y+g(y_0) $  for $0\leq y\leq y_0$. Then $f$ is continuous and decreasing.

From \eqref{eq:3.4}, for all $\epsilon>0$ there exist an $n(\epsilon)$ such that for  $\max\{y_0,n(\epsilon)\} <n\leq y<n+1$ we have $\mathbb{P}(\mathbf{Y} \leq y) = \mathbb{P}(\mathbf{Y} \leq n)$ and
\begin{equation}
1-((1-\underline{\beta}\rho)+\epsilon)f(y-1) \leq 
\mathbb{P}(\mathbf{Y} \leq y) \leq 1-((1-\overline{\beta}\rho)-\epsilon)f(y)
\label{eq:3.9b}
\end{equation}
For $\overline{\beta}\geq 0$, insertion of 
\begin{equation}
y=f^{-1}({\rm e}^{-x} / [(1-\underline{\beta}\rho)k]), \quad y=f^{-1}({\rm e}^{-x} / [(1-\overline{\beta}\rho)k]) ,
\label{eq:3.10}
\end{equation}
in the left-hand and right-hand equalities in \eqref{eq:3.9b} yields \eqref{eq:thm4-1} and \eqref{eq:thm4-2}.

If $\underline{\beta}>0$ then for all $\underline{\beta}\rho>\epsilon>0$ there exists an $y(\epsilon)$ such that $(\underline{\beta}-\epsilon)g(y)<g(y+1)<(\overline{\beta}-\epsilon)g(y)$, which is inherited by $f(y)$ for $y\geq\max\{y_0,y(\epsilon)\}$. For $\max\{y_0,n(\epsilon),y(\epsilon)\} < n \leq y<n+1$ we now obtain 
\[
1-((1-\underline{\beta}\rho)+\epsilon)/(1-\underline{\beta}\rho) f(y-1) \leq 
\mathbb{P}(\mathbf{Y} \leq y) \leq 1-((1-\overline{\beta}\rho)-\epsilon)f(y). 
\]
Insertion of $y$ as specified in \eqref{eq:3.10} yields \eqref{eq:thm4-2} and \eqref{eq:thm4-3}.  \hfill $\Box$ \\

If $\beta$ exists, and $0 < {\beta} \rho <1$, then there exist constants  $\{a_k>0,~b_k\}_{k\geq 1}$, $a_k$ bounded, $b_k\rightarrow \infty$ as $k\rightarrow \infty$, such that for  $-\infty < x < \infty$
\begin{equation}
{\rm e}^{-({\beta}\rho)^{-1}{\rm e}^{-x}} \leq \liminf_{k \rightarrow \infty} \mathbb{P}(\mathbf{Y}^{(k)} \leq a_k x + b_k) \leq
\limsup_{k \rightarrow \infty} \mathbb{P}(\mathbf{Y}^{(k)} \leq a_k x + b_k) \leq {\rm e}^{-{\rm e}^{-x}}, 
\label{eq:3.12}
\end{equation}
see \cite{Anderson1970}, with constants  such 
that 
\begin{equation}
\lim_{k\rightarrow \infty} f^{-1} ({\rm e}^{-y}/k) / (a_ky+b_k)=1.
\label{eq:4.1}
\end{equation}
This result  \eqref{eq:3.12} is not valid for $\overline{\beta}\neq\underline{\beta}$, as the arguments of $f^{-1}$ in the $\liminf$ and $\limsup$ differ in Theorem \ref{Thm:3.5}. In this case, constants  $\{a_k>0,~b_k\}_{k\geq 1}$ for the $\liminf$ and $\limsup$ need not coincide. 

The result \eqref{eq:3.12} cannot be strengthened, that is $\lim_{k \rightarrow \infty} \mathbb{P}(\mathbf{Y}^{(k)} \leq a_k x + b_k)$ does not exist, see \cite[Theorem 2.3]{Serfozo1988}. We may show that every sequence of sample maxima has a subsequence with limit between the upper and lower bounds in \eqref{eq:3.12} for fixed $x$. Stochastic compactness makes this precise and shows that this limit yields a non-degenerate distribution. 
\begin{Definition}[Stochastically compact (\cite{DeHaanResnick1984})]
The sequence of sample maxima $\mathbf{Y}^{(k)}$ is stochastically compact if there exist constants  $\{a_k>0,~b_k\}_{k\geq 1}$ such that every sequence $\{(\mathbf{Y}^{n(k)} - b_{n(k)})/a_{n(k)}\}_{k\geq 1}$ contains a subsequence whose distributions converge weakly to a non-degenerate probability distribution. 
Such a limit distribution is called a partial limit distribution.
We will also call the distribution of $\mathbf{Y}$ stochastically compact if the above holds. 
\end{Definition}
\begin{Theorem}\label{Thm:3.7}
If $\beta$ exists, and $0<\beta\rho <1$, then  $\mathbb{P}(\mathbf{Y} \leq y)$ is stochastically compact. The possible partial limit distributions are $G(x)=\exp (-\exp(-x+\epsilon))$, $-\infty<x<\infty$, with $\log(\beta\rho ) \leq \epsilon \leq 0$.
\end{Theorem}
 {\bf Proof} We will show that $\mathbb{P}(\mathbf{Y} \leq y)$ satisfies the sufficient conditions in \cite[Theorem 4]{DeHaanResnick1984}: For some $\delta >1$,
 \begin{align}
 & \int_x^{\infty} (1-\mathbb{P}(\mathbf{Y} \leq y))^{\delta-1}dy < \infty, \label{eq:3.14}\\
 & \liminf_{x \rightarrow \infty} \frac{\int_x^{\infty} (1-\mathbb{P}(\mathbf{Y} \leq y))^{\delta}dy}{ (1-\mathbb{P}(\mathbf{Y} \leq x))\int_x^{\infty} (1-\mathbb{P}(\mathbf{Y} \leq y))^{\delta-1}dy}>0, \label{eq:3.15} \\
 & \limsup_{x \rightarrow \infty} \frac{\int_x^{\infty} (1-\mathbb{P}(\mathbf{Y} \leq y))^{\delta}dy}{ (1-\mathbb{P}(\mathbf{Y} \leq x))\int_x^{\infty} (1-\mathbb{P}(\mathbf{Y} \leq y))^{\delta-1}dy}<1. \label{eq:3.16}
 \end{align}
These conditions are satisfied as a consequence of  \eqref{eq:3.5}. Let $0< \epsilon < \min\{\beta\rho, (1-\beta\rho)\}$. From  \eqref{eq:3.5}  there exists an $n(\epsilon)$ such that for all $ n>n(\epsilon)$ we have 
$((1-\beta\rho)-\epsilon)\psi(n)\rho^n \leq 
1-\mathbb{P}(\mathbf{Y} \leq n) \leq ((1-\beta\rho)+\epsilon)\psi(n)\rho^n
$,
and 
$
(\beta\rho)-\epsilon)\psi(n)\rho^n < \psi(n+1)\rho^{n+1}< (\beta\rho)+\epsilon)\psi(n)\rho^n
$. Combination of these equalities and observing that $\mathbb{P}(\mathbf{Y} \leq y) = \mathbb{P}(\mathbf{Y} \leq n)$ for $n\leq y < n+1$  shows that $\mathbb{P}(\mathbf{Y} \leq y)$ satisfies \eqref{eq:3.14} -- \eqref{eq:3.16}.
\hfill $\Box$ \\
 
 If $\beta$ exists, and $\beta=0$, from Lemma \ref{Lem:1} we obtain that
 \[
 \lim_{n\rightarrow\infty} \frac{1-\mathbb{P}(\mathbf{Y} \leq n-1)}{1-\mathbb{P}(\mathbf{Y} \leq n)} = \frac{1}{\rho}\lim_{n \rightarrow \infty} \frac{\psi(n-1)}{\psi(n)} = \infty, 
 \]
and \cite[Corollary 4]{DeHaanResnick1984} implies that $\mathbb{P}(\mathbf{Y} \leq y)$ is not stochastically compact. In this case, the bounds on the limiting distribution in Theorem \ref{Thm:3.5} do not enable characterisation of the limiting behaviour of $\mathbb{P}(\mathbf{Y}^{(k)} \leq a_k x + b_k) $.

The relation between \eqref{eq:Rem34-1} and \eqref{eq:Rem34-2}, suggests that  if $\beta$ exists, and $ {\beta} \rho > 1$  limiting results  for the limiting distribution of  $\mathbf{Y}^{(k)}$ conditioned on $\mathbf{Y}^{(k)}< \infty$ are similar 
to those of Theorem \ref{Thm:3.7}. 
 \begin{Theorem}\label{Thm:3.9}
If $\beta$ exists, and $1<\beta\rho <\infty$, then  $\mathbb{P}(\mathbf{Y} \leq y|\mathbf{Y} \leq \infty)$ is stochastically compact. The possible partial limit distributions are $G(x)=\exp (-\exp(-x+\epsilon))$, $-\infty<x<\infty$, with $-\log(\beta\rho ) \leq \epsilon \leq 0$.
\end{Theorem}
 {\bf Proof}
Follows the lines of the proof of Theorem \ref{Thm:3.7}, but now using the inequalities
$
(\frac{1}{(\beta\rho-1)B_*(B_*-1)}-\epsilon) \frac{1}{\psi(n)\rho^n} \leq 1- \mathbb{P}(\mathbf{Y} \leq y|\mathbf{Y} \leq \infty) \leq (\frac{1}{(\beta\rho-1)B_*(B_*-1)}+\epsilon) \frac{1}{\psi(n)\rho^n}
$
and
$(\frac{1}{\beta\rho}-\epsilon)\frac{1}{\psi(n)\rho^n} < \frac{1}{\psi(n+1)\rho^{n+1}}
< (\frac{1}{\beta\rho}+\epsilon)\frac{1}{\psi(n)\rho^n}$. \hfill $\Box$ \\

If  $\beta$ exists, and $\beta=\infty$, we obtain from \cite[Corollary 4]{DeHaanResnick1984} that $\mathbb{P}(\mathbf{Y} \leq y|\mathbf{Y} \leq \infty)$ is not stochastically compact.

\subsection{Almost sure convergence of sample maxima}
We will assume that constants  $\{a_k>0,~b_k\}_{k\geq 1}$, $a_k$ bounded, $b_k\rightarrow \infty$ as $k\rightarrow \infty$ exist such that \eqref{eq:4.1} is satisfied. Observe that  $\{b_k\}_{k\geq 1}$ may then be chosen as $b_k=f^{-1}(1/k)$. 

The following result is a corollary of Theorem \ref{Thm:3.5}. 
\begin{Corollary}\label{Cor:4.1} Assume that $0 \leq \overline{\beta}\rho <1$, and that the constants  $\{a_k>0,~b_k\}_{k\geq 1}$ satisfy \eqref{eq:4.1}. Then $\mathbf{Y}^{(k)} / b_k$ converges for $k\rightarrow \infty$ in probability to 1.
\end{Corollary}
The proof is omitted as it  follows standard lines as presented in \cite[Section 4.1]{Galambos1987}, but deviates as  the arguments $f^{-1}({\rm e}^{-x} / [(1-\underline{\beta}\rho)k]+1)$ in \eqref{eq:thm4-1}  and $f^{-1}({\rm e}^{-x} / [(1-\overline{\beta}\rho)k])$ in \eqref{eq:thm4-2} are not identical. 

If $\beta$ exists, and $0<\beta\rho<1$, the result of Corollary \ref{Cor:4.1} may be strengthened. The proof follows the lines of the proof of \cite[Theorem 4.4.4]{Galambos1987}, with the adaptation to incorporate the difference in the constants $b_k$. 
\begin{Theorem}\label{Thm:4.2} Assume that $\beta$ exists, and $0<\beta\rho<1$. Let $b_k=\log k / \log(1/\beta\rho)$. Then
\[
\mathbb{P}(\lim_{k\rightarrow\infty} \mathbf{Y}^{(k)} / b_k =1)=1.
\]
\end{Theorem}
{\bf Proof} Observe that $\sum_{n=1}^{\infty} (1-\mathbb{P}(\mathbf{Y} \leq tb_k)$ diverges for all $t\leq 1$ and converges for all $t>1$, so that \cite[Theorem 4.4.1]{Galambos1987} implies that 
$\mathbb{P}(\limsup_{k\rightarrow\infty} \mathbf{Y}^{(k)} / b_k =1)=1$. 

Observe that $\sum_{n=1}^{\infty}(1-\mathbb{P}(\mathbf{Y} \leq tb_k)) \exp(-k(1-\mathbb{P}(\mathbf{Y} \leq tb_k)) )$ converges for all $t <  1$, so that \cite[Theorem 4.3.3]{Galambos1987} with $u_k=tb_b$ implies that $\mathbb{P}(\liminf_{k\rightarrow\infty} \mathbf{Y}^{(k)} / b_k =1)=1$. \hfill $\Box$ \\

If $0<\underline{\beta}<\overline{\beta}$ we have  the bounds $\mathbb{P}(\limsup_{k\rightarrow\infty} \mathbf{Y}^{(k)} / \underline{b}_k \geq 1)=\mathbb{P}(\limsup_{k\rightarrow\infty} \mathbf{Y}^{(k)} / \overline{b}_k \leq 1)=1$, with $\underline{b}_k=\log k / \log(1/\underline{\beta}\rho) <   \log k / \log(1/\overline{\beta}\rho)=\overline{b}_k$, which is insufficient to conclude the first part of the proof above. If $\overline{\beta}=0$ both parts of the proof above fail. In that case we require the specific form of $f$ to obtain convergence results for $\mathbf{Y}^{(k)}$. Below, we consider the special case $f(x) \sim (\rho/x)^x$ as a natural extension of the geometric case $f(x) \sim \rho^x$. This result includes the infinite-server queue. 
\begin{Theorem}\label{Thm:4.3}
Let $f(x)=\alpha \gamma^x(1/x)^{x+\delta}$, $x>0$, with $\alpha>0$, $\gamma>0$. Let $b_k=f^{-1}(1/k)$.  If $\lim_{n\rightarrow\infty}\mathbb{P}( \mathbf{Y}\leq n) / f(n)=1$, then 
\[
\mathbb{P}(\lim_{k\rightarrow\infty} \mathbf{Y}^{(k)} / b_k =1)=1.
\]
\end{Theorem}
{\bf Proof} Exploring the form of $f$ and $b_k$, we readily find that 
$\sum_{n=1}^{\infty} (1-\mathbb{P}(\mathbf{Y} \leq tb_k)$ diverges for all $t\leq 1$ and converges for all $t>1$, and that $\sum_{n=1}^{\infty}(1-\mathbb{P}(\mathbf{Y} \leq tb_k)) \exp(-k(1-\mathbb{P}(\mathbf{Y} \leq tb_k)) )$ converges for all $t <  1$, so that we may apply \cite[Theorem 4.4.1]{Galambos1987} and \cite[Theorem 4.3.3]{Galambos1987}. \hfill $\Box$

\section{Examples}\label{sec:4}
\subsection{The multi-server queue}
The birth-death process with birth rates $\lambda(n)=\lambda$, $n=0,1,2,\ldots$, and death rates $\mu(n)=\min(n,s)\mu$, $n=1,2,\ldots$,  records the evolution of the number of customers in the  $M(\lambda)/M(\mu)/s$-queue, the multi-server queue with $s$ servers, that has Poisson arrivals with rate $\lambda$ and exponential service requirements with rate $\mu$. For $n=0,1,\ldots$,
\[
\psi(n)=\phi(n)=\left\{ \begin{array}{ll} \frac{1}{n!}, & n \leq s, \\[2mm]
\frac{s^{s-n}}{s!} , & n>s, 
\end{array}\right.
\]
If $\rho/s<1$, we have $\pi_{\phi}(n)=\pi_{\psi}(n)=B_{\phi} \phi(n)\rho^n$, $n=0,1,2,\ldots$. Clearly, $\beta$ exists, and $\beta = 1/s$.

Lemmas \ref{Lem:1} -- \ref{Lem:3} show that 
\begin{align*}
&  \lim_{n\rightarrow \infty} \left({\frac{s}{\rho}}\right)^n[1-\mathbb{P}(\mathbf{Y} \leq n)]= \frac{s^s}{s!}(1-\rho/s)  && \rho/s <1 , \\
 &\lim_{n\rightarrow \infty} n [1-\mathbb{P}(\mathbf{Y} \leq n)] =   \frac{s^s}{s!} & & \rho / s =1, \\
& \lim_{n\rightarrow \infty}  \left(\frac{\rho}{s}\right)^n[1-\mathbb{P}(\mathbf{Y} \leq n | \mathbf{Y}<\infty)] = [\frac{s^s}{s!}(\frac{\rho}{s}-1)  B_*(B_*-1)]^{-1},  &  & \rho / s >1.
\end{align*}
If $\beta \rho = \rho/s <1$, Theorem \ref{Thm:3.7} shows that $\mathbb{P}(\mathbf{Y} \leq y)$ is stochastically compact. The constants are $a_k = 1/ \log(s/\rho)$, $b_k=(\log k  )/ \log(s/\rho)$. The partial limit distribution is $G(x)=\exp (-\exp(-x+\epsilon))$, $-\infty<x<\infty$, with $\log(\rho / s) \leq \epsilon \leq 0$. Theorem \ref{Thm:4.2} shows that $\mathbf{Y} ^{(k)}/\log k$ converges to $\log (s/\rho)$ almost surely. 


\subsection{The infinite-server queue}\label{sec:4.2}
The birth-death process with birth rates $\lambda(n)=\lambda$, $n=0,1,2,\ldots$, and death rates $\mu(n)=n\mu$, $n=1,2,\ldots$,  records the evolution of the number of customers in the  $M(\lambda)/M(\mu)/\infty$-queue, the infinite-server queue, that has Poisson arrivals with rate $\lambda$ and exponential service requirements with rate $\mu$. For $n=0,1,\ldots$,
\[
\psi(n)=\phi(n)= \frac{1}{n!}.
\]
We have $\pi_{\phi}(n)=\pi_{\psi}(n)=B_{\phi} \phi(n)\rho^n$, $n=0,1,2,\ldots$. Clearly, $\beta$ exists, and $\beta = 0$.

Lemma \ref{Lem:1} shows that 
\[
 \lim_{n\rightarrow \infty} \frac{n!}{\rho^n} [1-\mathbb{P}(\mathbf{Y} \leq n)]= 1.
\]
This result does not appear as the limiting results for $s\rightarrow \infty$ in the multi-server queue. The factor $n!$ complicates the analysis of the limiting behaviour.  $\mathbb{P}(\mathbf{Y} \leq y)$ is not stochastically compact. Bounds on the limiting distribution in Theorem \ref{Thm:3.5} do not enable characterisation of the limiting behaviour of $\mathbb{P}(\mathbf{Y}^{(k)} \leq a_k x + b_k) $. 

Stirling's formula for $n!$ shows that the function $f(x)$ in  \eqref{eq:4.1} is  $f(x)=\frac{1}{\sqrt{2 \pi}} (\rho e)^x (1/x)^{x+{1\over 2}}$. The constants are $b_k =\log k / \log\log k$, and  $a_k=(\log b_k +1/(2b_k)-\log \rho)^{-1}$. Theorem \ref{Thm:4.3} shows that $\mathbf{Y} ^{(k)}/(\log k / \log\log k)$ converges to 1 almost surely.

\subsection{Kelly-Whittle network of queues}
The Markov chain $(X(t),~t\geq 0)$ at state space $S= \{0,1,2,\ldots\}^J$, with states ${\bf n}=(n_1,\ldots,n_J)$, and transition rates, for ${\bf n}' \neq {\bf n}$,
\[
q({\bf n},{\bf n}') = \left\{ \begin{array}{ll} {\displaystyle \mu_i
\frac{\psi({\bf n}-{\bf e}_i)}{\phi({\bf n})}p_{ij}}, & \mbox{if~}  {\bf n}' = {\bf n}-{\bf e}_i+{\bf e}_j, ~ i,j=0,\ldots,J, \\ [3mm]
0, & \mbox{otherwise,}
\end{array}
\right. 
\label{ch4:trKW}
\]
where ${\bf e}_i$ is the $i$-th unit vector, $i=1,\ldots,J$, ${\bf e}_0$ the zero-vector, $\psi:S\rightarrow (0,\infty)$ and $\phi:S\rightarrow (0,\infty)$, records the evolution of the number of customers in the queues of an open {Kelly-Whittle network} (\cite{BoucherievanDijk1991,HendersonTaylor1990}). 
If the routing matrix $P=(p_{ij},~i,j=0,\ldots,J)$  is irreducible,  the traffic equations $
\lambda_j= p_{0j} + \sum_{i=1}^{J} \lambda_i p_{ij}$, $j=1,\ldots,J
$, have a unique solution. Let $\rho_j= \lambda_j / \mu_j$, $j=1,\ldots,J$. If
$
B_{\phi,KW}^{-1} = \sum_{{\bf n} \in S} \phi({\bf n}) \prod_{j=1}^J (\mu_0\rho_j)^{n_j} < \infty,
$
then 
$(X(t),~t\geq 0)$ is positive recurrent with unique limiting distribution
$
\pi_{\phi}({\bf n})= B_{\phi,KW} \phi({\bf n}) \prod_{j=1}^J (\mu_0\rho_j)^{n_j}$, $ {\bf n} \in S
$, 
see \cite{BoucherieNorton}.
If
$
B_{\psi,KW}^{-1} = \sum_{{\bf n} \in S} \psi({\bf n}) \prod_{j=1}^J (\mu_0\rho_j)^{n_j} < \infty,
$
then $(X(t),~t\geq 0)$ is recurrent with limiting distribution at arrival epochs of a customer at a station $\pi_{\psi}({\bf n})= B_{\psi,KW}\psi({\bf n}) \prod_{j=1}^J (\mu_0\rho_j)^{n_j}$, $ {\bf n} \in S$, see \cite{Boucheriearrival}.

We are interested in the generic cycle maximum $\mathbf{Y} $ of the total number of customers in the network, $N=\sum_{i=1}^Jn_i$. 
To this end, for $N=0,1,2,\ldots$, let 
\[ 
\Psi(N) =
  \sum_{\{{\bf n} : \sum_{i=1}^J n_i =N\}} \psi({\bf n}) \prod_{j=1}^J \rho_j^{n_j} ,\quad
\Phi(N)=
  \sum_{\{{\bf n} : \sum_{i=1}^J n_i =N\}} \phi({\bf n} ) \prod_{j=1}^J \rho_j^{n_j}  \label{defPhi}.
\]
Norton's theorem shows that the birth-death process at $S=\{0,1,2,\ldots\}$ with birth and death rates 
\[
\lambda(N)=\mu_0 \frac{\Psi(N)}{\Phi(N)}, \mbox{ and }  \mu(N)=\frac{\Psi(N-1)}{\Phi(N)}
\]
 records the evolution of the total number of customers in the network of queues, see \cite[Theorem 6.4]{BoucherieNorton}. Thus, the generic cycle maximum $\mathbf Y$ has distribution $\mathbb{P}(\mathbf{Y} \leq N)= 1-\left[ \sum_{i=0}^N (\Psi(i)\mu_0^i)^{-1}\right]^{-1}$, $N=1,2,\ldots$. 
 
 The tail behaviour of the cycle maximum is determined by $\Psi(N)=  \sum_{\{{\bf n} : \sum_{i=1}^J n_i =N\}} \psi({\bf n}) \prod_{j=1}^J \rho_j^{n_j}$, which coincides with the normalising constant of a closed Kelly-Whittle network at state space $S_N=\{{\bf n} : \sum_{i=1}^J n_i =N\}$, with transition rates 
$q({\bf n},{\bf n}') =
\mu_i
\frac{\psi({\bf n}-{\bf e}_i)}{\psi({\bf n})}p_{ij}$, 
if 
${\bf n}' = {\bf n}-{\bf e}_i+{\bf e}_j$, $ i,j=1,\ldots,J$.
The limiting distribution is $\pi_{\psi,N}({\bf n}) = \Psi(N)^{-1}\psi({\bf n}) \prod_{j=1}^J \rho_j^{n_j} $. The  normalising constant for closed networks is well-studied in literature. Its asymptotic behaviour is available for several special cases, including networks comprised of multi-server and infinite-server queues. 

\subsubsection{Network of infinite-server queues}
If all queues are infinite-server queues, then $\phi({\bf n})=\psi({\bf n}) = \prod_{i=1}^J 1 / n_i!$, and $\Phi(N)=\Psi(N)=(
\sum_j \rho_J)^N / N!$. 
Alternatively, we may write the birth and death rates as $
\lambda(N)=\lambda \frac{\widehat{\Psi}(N)}{\widehat{\Phi}(N)}$,  and $\mu(N)=\mu\frac{\widehat{\Psi}(N-1)}{\widehat{\Phi}(N)}$, with $\lambda=\mu_0$, $\mu=1 / \sum_j \rho_j$, $\widehat{\Phi}(N)=\widehat{\Psi}(N)=1/N!$. Section \ref{sec:4.2} presents the limiting results. 

\subsubsection{Network of single-server queues}
If all queues are single-server queues, then $\phi({\bf n})=\psi({\bf n}) = 1$, and $\Phi(N)=\Psi(N)= 
  \sum_{\{{\bf n} : \sum_{i=1}^J n_i =N\}}  \prod_{j=1}^J \rho_j^{n_j} $. If all $\rho_j$ are distinct, from \cite[Theorem 1]{Harrison1985}, 
  $
  \Psi(N)= \sum_{j=1}^J \rho_j^{N+J-1} / \prod_{i \neq j} (\rho_j-\rho_i)
  $.

  We obtain that $\beta$ exists, and $\beta=\max\{\rho_1,\ldots,\rho_J\}$. If $\beta\mu_0<1$,  then $\lim_{N\rightarrow\infty}(\Psi(N) \mu_0^N)/ (\beta\mu_0)^N = c$, for some constant $c$. Theorem \ref{Thm:3.7} shows that $\mathbb{P}(\mathbf{Y} \leq y)$ is stochastically compact. The constants are $a_k=1/\log(1/(\beta\mu_0))$ and $b_k=\log k / \log(1/(\beta\mu_0))$.
 The partial limit distribution is $G(x)=\exp (-\exp(-x+\epsilon))$, $-\infty<x<\infty$, with $\log(\beta\mu_0) \leq \epsilon \leq 0$. Theorem \ref{Thm:4.2} shows that $\mathbf{Y} ^{(k)}/\log k$ converges to $\log (1/(\beta\mu_0))$ almost surely.

\subsubsection{Network of multi-server and infinite-server queues}
If the network contains single-server, multi-server and infinite-server queues, an explicit expression for $\Phi(N)=\Psi(N)$ is not available in closed form. Asymptotic results are available in literature, see e.g.,
\cite[Theorem 4]{GeorgeXiaSquillante2012}, where it shown that $\lim_{N\rightarrow \infty } \Psi(N) / \max\{\nu_1,\ldots,\nu_J\}^N)=c$, for some constant $c$, with $\nu_i= \rho_i / s_i$ all distinct, and $s_i$ the number of available servers at the station. Note that the term $\max\{\nu_1,\ldots,\nu_J\}^N$ does not appear in \cite[Theorem 4]{GeorgeXiaSquillante2012} as in that theorem the $\nu_j$ are scaled such that $\max\{\nu_1,\ldots,\nu_J\}=1$.  

We obtain that $\beta$ exists, and $\beta=\max\{\nu_1,\ldots,\nu_J\}$. If $\beta\mu_0<1$,  then $\lim_{N\rightarrow\infty}(\Psi(N) \mu_0^N)/ (\beta\mu_0)^N = c$, for some constant $c$.
Theorem \ref{Thm:3.7} shows that $\mathbb{P}(\mathbf{Y} \leq y)$ is stochastically compact. The constants are $a_k=1/\log(1/(\beta\mu_0))$ and $b_k=\log k / \log(1/(\beta\mu_0))$. The partial limit distribution is $G(x)=\exp (-\exp(-x+\epsilon))$, $-\infty<x<\infty$, with $\log(\beta \mu_0) \leq \epsilon \leq 0$. Theorem \ref{Thm:4.2} shows that $\mathbf{Y} ^{(k)}/\log k$ converges to $\log (1/(\beta\mu_0))$ almost surely. 

  \subsubsection{Network with identical load parameters}
We may also evaluate $\Psi(N)$ if some of the load parameters coincide. For the network of single server queues explicit expressions are available in \cite[Theorem 2]{Harrison1985}. For the network containing single-server, multi-server and infinite-server stations an asymptotic expression is available in \cite[Theorem 4]{GeorgeXiaSquillante2012}, where it shown that $\lim_{N\rightarrow \infty } \Psi(N) / (N^{|B|-1}\max\{\nu_1,\ldots,\nu_J\}^N)=c$, for some constant $c$, with $B = \{j : \nu_j = \max\{\nu_1,\ldots,\nu_J\}\}$. If two parameters coincide, we may use the Lambert $W$ function (see \cite{CorlessGonnetHare3JeffreyKnuth1996}) to approximate the inverse of $f(y)=y (\beta\mu_0)^y$, with $\beta = \max\{\nu_1,\ldots,\nu_J\}$. If $B$ parameters coincide the inverse of $f(y)=y^{|B|-1} (\beta\mu_0)^y$ must be evaluated.

\section{Concluding remarks}
The result presented in this paper build upon classical results from extremal value theory for discrete distributions (\cite{Anderson1970,Vervaat1973}), and birth-death processes (\cite{Asmussen1998,Cohen1982}). The results are applied to study convergence of the sequence of partial maxima for single queues and open Kelly-Whittle networks of queues. The latter results rely on an amenable expression for the normalising constant for suitably defined closed networks. Such expressions are available for networks of single-server, multi-server and infinite-server queues with distinct load parameters.

The results may readily be generalised to networks with state-dependent routing (\cite{BoucherieNorton}). For positive recurrent networks of quasi-reversible queues the birth and death rates may also be obtained using Norton's theorem, see \cite{BoucherieNortonquasireversible}.  Careful analysis of these birth and death rates is required to 
determine the asymptotic behaviour of these rates, and hence the limiting behaviour of the sequence of sample maxima.

\section*{Acknowledgements} The author is grateful to prof. J.W. Cohen for inspiring discussions on the dual relation between 
the discrete failure rate $
\mathbb{P}(\mathbf{Y} = n | \mathbf{Y} \geq n)$ and the blocking probability $
\mathbb{P}(\mathbf{X} = n | \mathbf{X} \leq n)$ that led to this research. The results for birth-death processes presented in the unpublished report (\cite{Boucheriemaximum}), are now extended to open networks of queues.

\bibliographystyle{unsrtnat}
\bibliography{references}
\end{document}